\documentclass[11pt]{amsart}
\usepackage{amsmath,amsfonts,amssymb,amsthm}
\usepackage[alphabetic]{amsrefs}
\usepackage{hyperref}
\usepackage{graphicx,color}
\usepackage{enumerate}
\usepackage{tikz}

\newtheorem{theorem}{Theorem}[section]
\newtheorem{lemma}[theorem]{Lemma}
\newtheorem{prop}[theorem]{Proposition}
\newtheorem{cor}[theorem]{Corollary}

\newtheorem{main}{Theorem}

\newtheorem{cmain}[main]{Corollary}
\newtheorem{qmain}[main]{Question}
\newtheorem{conjmain}[main]{Conjecture}

\theoremstyle{definition}

\newtheorem*{question}{Questions}

\newcommand{\N}{\mathbb{N}}
\newcommand{\Z}{\mathbb{Z}}
\newcommand{\R}{\mathbb{R}}

\newcommand{\mk}{\medskip}
\newcommand{\bthm}{\begin{theorem}}
\newcommand{\ethm}{\end{theorem}}
\newcommand{\bmthm}{\begin{main}}
\newcommand{\emthm}{\end{main}}
\newcommand{\ben}{\begin{enumerate}}
\newcommand{\een}{\end{enumerate}}
\newcommand{\bit}{\begin{itemize}}
\newcommand{\eit}{\end{itemize}}
\newcommand{\f}{\frac}
\newcommand{\bp}{\begin{proof} \upshape}
\newcommand{\ep}{\end{proof}}
\newcommand{\beq}{\begin{eqnarray*}}
\newcommand{\eeq}{\end{eqnarray*}}
\newcommand{\ra}{\rightarrow}
\newcommand{\blem}{\begin{lemma}}
\newcommand{\elem}{\end{lemma}}
\newcommand{\bpro}{\begin{prop}}
\newcommand{\epro}{\end{prop}}
\newcommand{\bs}{\symbol{92}}
\newcommand{\bcor}{\begin{cor}}
\newcommand{\ecor}{\end{cor}}

\newcommand{\bmcor}{\begin{cmain}}
\newcommand{\emcor}{\end{cmain}}
\newcommand{\bq}{\begin{question}}
\newcommand{\eq}{\end{question}}
\newcommand{\bmq}{\begin{qmain}}
\newcommand{\emq}{\end{qmain}}
\newcommand{\bmconj}{\begin{conjmain}}
\newcommand{\emconj}{\end{conjmain}}
\newcommand{\eps}{\varepsilon}

\newcommand{\st}{\, | \,}

\begin{document}
	
	\title{Automorphisms and subdivisions of Helly graphs}
	
	
	\author[T.~Haettel]{Thomas Haettel}

\let\thefootnote\relax\footnotetext{Thomas Haettel, thomas.haettel@umontpellier.fr, IMAG, Univ Montpellier, CNRS, France, and IRL 3457, CRM-CNRS, Universit\'{e} de Montr\'{e}al, Canada.}

	\date{\today}

	\begin{abstract}

		
		We study Helly graphs of finite combinatorial dimension, i.e. whose injective hull is finite-dimensional. We describe very simple fine simplicial subdivisions of the injective hull of a Helly graph, following work of Lang. We also give a very explicit simplicial model of the injective hull of a Helly graph, in terms of cliques which are intersections of balls.
		
		We use these subdivisions to prove that any automorphism of a Helly graph with finite combinatorial dimension is either elliptic or hyperbolic. Moreover, every such hyperbolic automorphism has an axis in an appropriate Helly subdivision, and its translation length is rational with uniformly bounded denominator.
	\end{abstract}
	
	\maketitle

\let\thefootnote\relax
	
	\section{Introduction}

A connected graph such that any family of pairwise intersecting balls has a non-empty global intersection is called a Helly graph. Such graphs appear to play an increasing role in geometric group theory, as many groups have interesting actions on Helly graphs, most notably Gromov-hyperbolic groups, cubulated groups, braid groups and some higher rank lattices (see notably~\cite{lang,helly_groups,hoda:crystallographic,osajda_valiunas,haettel_injective_buildings,haettel_hoda_petyt,haettel_helly_kpi1,haettel_osajda_locally_elliptic,haettel_huang_garside_artin_product_Z}).

\mk

One of the most natural questions, when studying a metric space which has some form of nonpositive curvature, is to study the possible individual isometries.
	
\mk

In order to study automorphisms of a Helly graph $X$, we are interested in finding a nice combinatorial structure on the injective hull $E(X)$. Such a description has been carried out by Lang in~\cite{lang}, and we present a slight modification of his construction, see Theorem~\ref{thm:orthoscheme_complex} for the precise statement. Recall that the \emph{combinatorial dimension} of $X$ is the dimension of its injective hull $E(X)$. Note that any group acts on the Helly hull of its Cayley graph, so one needs to restrict the class of Helly graphs we will be considering: we will hence mostly be considering Helly graphs with finite combinatorial dimension.

\bmthm[Orthoscheme complex of a Helly graph] \label{mthm:orthoscheme_complex}\
Let $X$ denote a Helly graph with finite combinatorial dimension. For each $N \geq 1$, there exists a simplicial structure on the injective hull $E(X)$ of $X$, denoted $O_N X$ and called the ($N^\text{th}$) \emph{orthoscheme subdivision complex} of $X$, satisfying the following:
\bit
\item Each simplex of $O_N X$ is isometric to the standard $\ell^\infty$ orthosimplex with edge lengths $\f{1}{2N!}$.
\item The vertex set $X'_N$ of $O_N X$, endowed with the induced distance, is a Helly graph (with edge lengths $\f{1}{2N!}$), containing isometrically $X$, called the ($N^\text{th}$) Helly subdivision of $X$. Moreover, we have
$$X'_N = \left\{p \in E(X) \st \forall x \in X, d(p,x) \in \f{1}{2N!}\N\right\}.$$
\eit
\emthm

\mk

We also obtain a very explicit description of the first subdivision. Let us recall that, in this article, a \emph{clique} of a graph is the vertex set of a complete subgraph. Moreover, let us say that a clique is \emph{round} if it is an intersection of balls.

\bmthm[First subdivision]
Let $X$ denote a Helly graph. The following graphs are naturally isomorphic:
\bit
\item The graph with vertex set
$$E(X) \cap \left( \f12 \N\right)^X,$$
with an edge between $f,g \in E(X)$ if and only if $d_\infty(f,g)=\f12$.
\item The graph with vertex set
$$\{\mbox{round cliques of }X\} = X \cup \{\mbox{non-empty intersections of maximal cliques of }X\},$$
with an edge between $\sigma,\tau \subset X$ if and only if $\sigma \cap \tau \neq \emptyset$ and $\sigma \cup \tau$ is a clique of $X$.
\eit
This graph coincides with the \emph{first Helly subdivision} $X'=X'_1$ of $X$ described in Theorem~\ref{mthm:orthoscheme_complex}, and the natural map $X \ra X'$ is a $2$-homothetic embedding.
\emthm

One nice consequence is a very simple characterization of the combinatorial dimension of a Helly graph.

\bmcor \label{mcor:identification_combinatorial_dimension}
Let $X$ denote a Helly graph. Then the combinatorial dimension of $X$ coincides with the length of the longest chain of round cliques of $X$.
\emcor

In particular, this bounds easily the combinatorial dimension of Helly graphs with bounded valence.

\bmcor
Any Helly graph of valence at most $N$ has combinatorial dimension at most $N-1$.
\emcor

Moreover, many locally infinite Helly graphs can also be shown to have finite combinatorial dimension: for instance, every tree has combinatorial dimension at most $1$.

\mk

We will use the orthoscheme subdivision to study automorphisms of Helly graphs. One key property of CAT(0) spaces is the classification of isometries into elliptic, parabolic and hyperbolic (see~\cite[Definition~6.3]{bridson_haefliger}).

\mk

In this article, we prove a similar classification for automorphisms of Helly graphs. We say that an automorphism (or a group of automorphisms) of a Helly graph is \emph{elliptic} if it stabilizes a clique of $X$. We say that an automorphism $g$ of a Helly graph is hyperbolic if the orbit map $n \in \Z \mapsto g^n \cdot x$ is a quasi-isometric embedding. We refer to Section~\ref{sec:classification} for other characterizations of elliptic and hyperbolic automorphisms, and to Theorem~\ref{thm:dichotomy} for the precise statement.

\bmthm[Classification of automorphisms of Helly graphs] \label{mthm:classification_isometries}\

Let $X$ denote a Helly graph with finite combinatorial dimension $N$. Then any automorphism of $X$ is either elliptic or hyperbolic.

\mk

More precisely, any elliptic automorphism of $X$ fixes a vertex in the $N^\text{th}$ Helly subdivision $X'_N$ of $X$.

Every hyperbolic automorphism $g$ of $X$ has a combinatorial axis in the $N^\text{th}$ Helly subdivision $X'_N$ of $X$, i.e. there exists a vertex $x \in X'_N$ such that $(g^n \cdot x)_{n \in \Z}$ is a geodesic in $X'_N$.

\mk

In addition, every hyperbolic automorphism of $X$ has rational translation length, with denominator bounded above by $2N$.\emthm

This is a direct generalization (in the finite-dimensional case) of a result of Haglund stating essentially that any automorphism of a CAT(0) cube complex either fixes a point or translates a combinatorial geodesic (see~\cite[Theorem~1.4]{haglund} for the precise statement).

\mk

This also generalizes a theorem of Gromov for translation lengths of hyperbolic elements in a Gromov-hyperbolic group (see~\cite[8.5.S]{gromov_hyperbolic_groups}). Since Garside groups are Helly according to~\cite{huang_osajda_helly}, this implies a direct analogue of~\cite{lee_lee_garside_translation} for a very closely related translation length. This has consequences in particular for decision problems, following~\cite{lee_lee_garside_translation}, since the conjugacy problem is solvable for Helly groups (see~\cite{helly_groups}). 

\bmcor
Let $G$ denote a Helly group. The following problems are solvable for $G$.
\bit
\item The power problem: given infinite order elements $g,h \in G$, find $n \geq 1$ such that $h^n=g$.
\item The power conjugacy problem: given infinite order elements $g,h \in G$, find $n \geq 1$ such that $h^n$ is conjugate to $g$.
\eit
\emcor

This result also has a direct consequence concerning distortion. Recall that an element $g$ of a finitely generated group $G$ with a word metric $|\cdot|_G$ is undistorted if there exists $C>0$ such that $\forall n \in \N, |g^n|_G \geq nC$.

\bmcor[No distortion in Helly graphs]\

Let $X$ denote a Helly graph with finite combinatorial dimension, let $G$ denote a finitely generated group of automorphisms of $X$ and assume that some element $g$ of $G$ is not elliptic. Then $g$ is hyperbolic in $X$, has infinite order and is undistorted in $G$. 

More precisely, $g$ is uniformly undistorted : $\exists C>0, \forall n \in \N, |g^n|_G \geq nC$. Explicitly, if some orbit map for the action of $G$ on $X$ is $K$-Lipschitz, one may choose $C=\f{1}{2NK}$. \emcor

If a finitely generated group $G$ acts properly by automorphisms on a Helly graph with finite combinatorial dimension, we therefore deduce that $G$ has \emph{uniformly undistorted infinite cyclic subgroups} as defined by Cornulier in~\cite[Definition~6.A.3]{cornulier_commensurated}. See also~\cite{abbott_hagen_petyt_zalloum} for a related statement.

This applies in particular to all discrete subgroups of semisimple Lie groups over non-Archimedean local fields of types $A$, $B$, $C$ or $D$, see~\cite{haettel_injective_buildings}.

\mk

In particular, we deduce an obstruction to the existence of some actions on Helly graphs.

\bmcor
Finitely generated groups with distorted elements do not act properly on a Helly graph with finite combinatorial dimension.
\emcor

This applies notably to nilpotent groups that are not virtually abelian, to non-uniform irreducible lattices in real semisimple Lie groups of higher rank and to Baumslag-Solitar groups. This generalizes, in the nilpotent case, the fact that every solvable subgroup of a Helly group is virtually abelian, see~\cite{valiunas_abelian_helly}.

\mk

Note that, on the other hand, any finitely generated group acts properly by automorphisms on a Helly graph, the Helly hull of any Cayley graph. In the case of a group with distorted elements, we deduce that the Helly hull of a Cayley graph has infinite combinatorial dimension.

\mk

The case of non-uniform lattices is drastically different from the uniform one. Indeed, uniform lattices in semisimple Lie groups over local fields have nice actions on Helly graphs (in the non-Archimedean case, see Theorem~\ref{thm:buildings_helly} and~\cite{haettel_injective_buildings} for details, and  \cite{haettel_helly_kpi1}) and on injective metric spaces (in the Archimedean case, see~\cite{haettel_injective_buildings} for details).

\mk

We also prove a result about fixed point sets for a pair of elliptic subgroups of a Helly graph with finite combinatorial dimension, which is used in~\cite{haettel_osajda_locally_elliptic} with Damian Osajda in our study of locally elliptic actions on Helly graphs.

\mk

\textbf{Organization of the article:} 
In Section~\ref{sec:Helly}, we review classical results about Helly graphs and injective metric spaces, mostly following work of Lang. In Section~\ref{sec:subdivisions}, we describe nice simplicial subdivisions of Lang's cell structure on the injective hull of a Helly graph. In Section~\ref{sec:first_subdivision}, we give a very explicit description of the first subdivision of a Helly graph, using round cliques. In Section~\ref{sec:classification}, we use these subdivisions to prove the classification result of automorphisms of Helly graphs. In Section~\ref{sec:fixed_point_pair}, we study fixed point set of pairs of elliptic subgroups.

\mk

\textbf{Acknowledgments:} 
We would like to thank Giuliano Basso, Anthony Genevois, Damian Osajda and Urs Lang for many interesting discussions on this work. We would also like to thank the referee for precise and insightful comments.

\mk

The author was partially supported by French projects ANR-16-CE40-0022-01 AGIRA and ANR-22-CE40-0004 GOFR.

\tableofcontents	

\section{Helly graphs} \label{sec:Helly}

A connected graph $X$ is called \emph{Helly} if any family of pairwise intersecting combinatorial balls of $X$ has a non-empty global intersection. We will consider $X$ as its vertex set, and we will endow $X$ with induced graph metric. We refer the reader to~\cite{helly_groups} for a presentation of Helly graphs and Helly groups.

	One may think of Helly graphs as a very nice class of nonpositively curved, combinatorially defined spaces. Surprisingly enough, many nonpositive curvature metric spaces and groups have a very close relationship to Helly graphs or their non-discrete counterpart, injective metric spaces.
	
	For instance, the thickening of any CAT(0) cube complex is a Helly graph (see~\cite{bandelt_vandevel_superextensions}, and also~\cite[Corollary~3.6]{hruskawise:packing}). Lang showed that the any Gromov hyperbolic group acts properly cocompactly on the Helly hull of any Cayley graph (see~\cite{lang,helly_groups}). Huang and Osajda proved that any weak Garside group and any Artin group of type FC has a proper and cocompact action on a Helly graph (see~\cite{huang_osajda_helly}, and also~\cite{haettel_helly_kpi1}). Osajda and Valiunas proved that any group that is hyperbolic relative to Helly groups is Helly (see~\cite{osajda_valiunas}). Haettel, Hoda and Petyt proved that any hierarchically hyperbolic group, and in particular any mapping class group of a surface, has a proper and cobounded action on an injective metric space, see~\cite{haettel_hoda_petyt}.

\mk

Concerning Euclidean buildings, recall the following statement.

\bthm[Hirai, Chalopin et al, Haettel] \label{thm:buildings_helly}
The thickening of any Euclidean building of type $\tilde{A}$ extended, $\tilde{B}$, $\tilde{C}$ or $\tilde{D}$ is Helly.
\ethm

Hirai, and Chalopin et al. proved the case of Euclidean buildings of type $\tilde{A}$ extended and $\tilde{C}$, see~\cite{hirai_uniform_modular} and \cite{chalopin_chepoi_hirai_osajda}. In~\cite{haettel_injective_buildings} and \cite{haettel_helly_kpi1}, Haettel proved the statement for all Euclidean buildings of type $\tilde{A}$ extended, $\tilde{B}$, $\tilde{C}$ or $\tilde{D}$. There is an analogous result for classical symmetric spaces, see~\cite{haettel_injective_buildings} for a precise statement.

\mk

Recall that a geodesic metric space is called \emph{injective} if any family of pairwise intersecting closed balls has a non-empty global intersection. We refer the reader to~\cite{lang} for a presentation of injective metric spaces, and also the following result of Isbell.

\bthm[\cite{isbell}]
Let $X$ denote a metric space. Then there exists an essentially unique minimal injective space $E(X)$ containing $X$, called the injective hull of $X$.
\ethm

In~\cite{lang}, Lang gives a very explicit description of the injective hull of a metric space $X$: let
$$\Delta(X)=\{f:X \ra \R \mbox{ $1$-Lipschitz}, \forall x,y \in X, f(x) +f(y) \geq d(x,y)\},$$
endowed with the sup metric. An element $f \in \Delta(X)$ is called \emph{extremal} if
$$\forall x \in X, f(x)= \sup_{y \in X} d(x,y)-f(y).$$
Then we can state Lang's result.

\bthm \cite[Theorem~3.3]{lang} \label{thm:explicit_description_injective_hull}
Let $X$ denote a metric space, then the space
$$ E(X)=\{f \in \Delta(X) \mbox{ $f$ is extremal}\},$$
with the isometric embedding $e:x \in X \mapsto d(x,\cdot) \in E(X)$, is the injective hull of $X$.
\ethm

\begin{figure}
\begin{center}
\begin{tikzpicture}
\def \p {0.05}
\def \op {0.5}
\def \gris {black!10}
\draw[fill] (0,0) circle (\p) node(0) {};
\draw[fill] (3,0) circle (\p) node(1) {};
\draw[fill] (3,3) circle (\p) node(2) {};
\draw[fill] (3,3) + (20:3) circle (\p) node(3) {};

\draw[black,fill opacity=\op,fill=\gris] (0.center) -- (1.center) -- (3.center) -- (0.center);
\draw[black,fill opacity=\op,fill=\gris] (0.center) -- (2.center) -- (3.center) -- (0.center);
\draw[black,fill opacity=\op,fill=\gris] (0.center) -- (1.center) -- (2.center) -- (0.center);
\draw[black,fill opacity=\op,fill=\gris] (3.center) -- (1.center) -- (2.center) -- (3.center);
\node at ([yshift=-0.5cm]0) {\bfseries $e_0=(0,0,0)$};
\node at ([yshift=-0.5cm]1) {\bfseries $e_1=(1,0,0)$};
\node at ([yshift=0.7cm]2) {\bfseries $e_2=(1,1,0)$};
\node at ([yshift=0.5cm]3) {\bfseries $e_3=(1,1,1)$};

\end{tikzpicture}
\end{center}
\caption{The standard orthosimplex of dimension $3$ with edge lengths $1$.}
\label{fig:orthosimplex}
\end{figure}
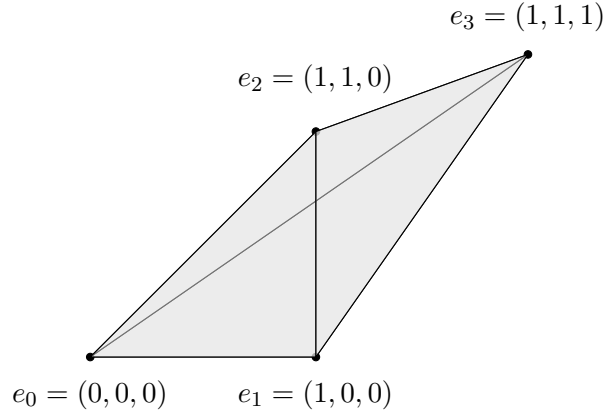

We will be mostly interested in the case where $X$ is the vertex set of a connected graph. Moreover, Lang describes a cell structure on the injective hull of a connected graph. We describe below a refinement of Lang's cell decomposition into orthosimplices. Recall that the \emph{standard orthosimplex} of dimension $n$ with edge lengths $\ell>0$ is the simplex of $\R^n$ with vertices $(0,\dots,0),(\ell,0,\dots,0),\dots,(\ell,\ell,\dots,\ell)$, see Figure~\ref{fig:orthosimplex}. We will endow this simplex with the standard $\ell^\infty$ metric on $\R^n$.

\mk

Recall that the \emph{combinatorial dimension} of a metric space $X$ is the dimension of its injective hull $E(X)$ (this has been defined by Dress, see~\cite{dress}). There are interesting examples of locally infinite Helly graphs with finite combinatorial dimension, such as thickenings of locally infinite, finite-dimen\-sional CAT(0) cube complexes.
%
%

\section{Helly subdivisions} \label{sec:subdivisions}

We now present a refinement of Lang's description of the cell structure on the injective hull of a connected graph (see~\cite{lang}).

\bthm \label{thm:orthoscheme_complex}
Let $X$ denote a Helly graph with finite combinatorial dimension. For each $N \geq 1$, there exists a simplicial structure on the injective hull $E(X)$ of $X$, denoted $O_N X$ and called the ($N^\text{th}$) \emph{orthoscheme subdivision complex} of $X$, satisfying the following:
\bit
\item Each simplex of $O_N X$ is isometric to the standard $\ell^\infty$ orthosimplex with edge lengths $\f{1}{2N!}$.
\item The vertex set $X'_N$ of $O_N X$, endowed with the induced distance, is a Helly graph (with edge lengths $\f{1}{2N!}$), containing isometrically $X$, called the ($N^\text{th}$) Helly subdivision of $X$. Moreover, we have
$$X'_N = \left\{p \in E(X) \st \forall x \in X, d(p,x) \in \f{1}{2N!}\N\right\}.$$
\item For any $p \in O_NX$ and for any simplex of $O_NX$ containing $p$ with vertices $x_1,\dots,x_n$ in $X'_N$, there exist unique  $t_1,\dots,t_n \geq 0$ such that\\ $t_1+\dots+t_n=1$ and
$$\forall x \in X, d(p,x)=\sum_{i=1}^n t_id(x_i,x).$$
We write $p=\sum_{i=1}^n t_ix_i$. More generally, if $p,p' \in O_NX$ are such that $p=\sum_{i=1}^n t_ix_i$ and $p'=\sum_{i'=1}^{n'} t'_{i'}x'_{i'}$, then
$$d(p,q) = \max_{x \in X} \sum_{i=1}^n \sum_{i'=1}^{n'} t_i t'_{i'} |d(x_i,x)-d(x'_{i'},x)|.$$
\eit
\ethm

Before passing to the proof, let us first explain why we want to consider $2N!$ and not $2^N$ for instance. Consider the Helly graph $\Gamma$ with vertex set $\Z^N$, with the standard Helly structure. Let $g$ denote the following automorphism of $\Gamma$:
$$g \cdot (x_1,x_2,\dots,x_N) = (x_2+1,x_3,x_4,\dots,x_N,x_1).$$
The automorphism $g$ is hyperbolic with translation length $\f{1}{N}$, hence if $N=3$ and $k$ is a power of $2$ for instance, then $g^k$ does not have a combinatorial axis in $\Gamma$.

\bp
According to~\cite[Theorem~4.5]{lang}, the injective hull $E(X)$ may be realized as an isometric subset of $\R^X$, and the injective hull $E(X)$ of $X$ has a natural cell decomposition satisfying the following. For each cell $C$ of $E(X)$, there is a finite set of vertices $x_1,\dots,x_n$ of $X$ such that the map
\beq C & \ra & \R^n \\
p & \mapsto & (d(p,x_1),\dots,d(p,x_n))\eeq
is an isometry (with the $\ell^\infty$ metric on $\R^n$) onto the compact convex subspace of $\R^n$ defined by inequalities of the type
$$ \pm d(\cdot,x_i) \pm d(\cdot,x_j) \leq D,$$
for some $1 \leq i < j \leq n$ and $D \in \Z$, and also of the type
$$\pm d(\cdot,x_i) \leq  D',$$
for some $1 \leq i \leq n$ and $D' \in \f{1}{2}\Z$.
In particular there is an affine structure on $C$. Moreover, for any $x \in X$, for any $p _1,\dots,p_k \in C$ and $t_1,\dots,t_k \geq 0$ such that $t_1+\dots+t_k=1$, we have
$$d(x,\sum_{i=1}^k t_ip_i) = \sum_{i=1}^k t_id(x,p_i).$$

\mk

\begin{figure}
\begin{center}
\begin{tikzpicture}[scale = 0.6]
\def \p {0.05}
\def \r {2}
\def \m {1}
\def \a {20}
\def \op {0.5}
\def \gris {black!10}
\draw[fill] (0,0) circle (\p) node(000) {};
\draw[fill] (4,0) circle (\p) node(100) {};
\draw[fill] (0,4) circle (\p) node(010) {};
\draw[fill] (4,4) circle (\p) node(110) {};
\draw[fill] (2,0) circle (\p) node(200) {};
\draw[fill] (0,2) circle (\p) node(020) {};
\draw[fill] (2,4) circle (\p) node(210) {};
\draw[fill] (4,2) circle (\p) node(120) {};
\draw[fill] (2,2) circle (\p) node(220) {};
\draw[fill] (2,2) + (\a:\m) circle (\p) node(222) {};
\draw[fill] (4,0) + (\a:\r) circle (\p) node(101) {};
\draw[fill] (0,4) + (\a:\r) circle (\p) node(011) {};
\draw[fill] (4,4) + (\a:\r) circle (\p) node(111) {};
\draw[fill] (4,0) + (\a:\m) circle (\p) node(102) {};
\draw[fill] (4,2) + (\a:\m) circle (\p) node(122) {};
\draw[fill] (4,2) + (\a:\r) circle (\p) node(121) {};
\draw[fill] (4,4) + (\a:\m) circle (\p) node(112) {};
\draw[fill] (0,4) + (\a:\m) circle (\p) node(012) {};
\draw[fill] (2,4) + (\a:\r) circle (\p) node(211) {};
\draw[fill] (2,4) + (\a:\m) circle (\p) node(212) {};

\draw (000.center) -- (100.center) -- (110.center) -- (010.center) -- (000.center);
\draw (000.center) -- (110.center);
\draw (100.center) -- (010.center);
\draw (200.center) -- (210.center);
\draw (020.center) -- (120.center);

\draw (100.center) -- (101.center) -- (111.center) -- (110.center);
\draw (100.center) -- (111.center) -- (010.center);
\draw (101.center) -- (110.center) -- (011.center);
\draw (010.center) -- (011.center) -- (111.center);
\draw (012.center) -- (112.center) -- (102.center);
\draw (210.center) -- (211.center);
\draw (120.center) -- (121.center);

\draw[blue, dashed] (000.center) -- (222.center);
\draw[blue,fill opacity=0.3,fill=blue] (220.center) -- (000.center) -- (200.center) -- (220.center) -- (222.center) -- (200.center);


\end{tikzpicture}
\end{center}
\caption{The partition of a cube in $\R^3$ into standard orthosimplices.}
\label{fig:partition}
\end{figure}
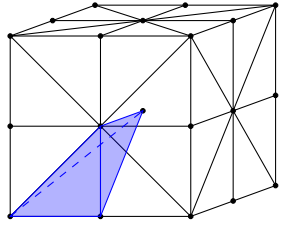

Note that the hyperplanes of $\R^n$
$$\left\{\pm x_i \pm x_j = D \st 1 \leq i < j, D \in \f{1}{N!}\Z\right\} \mbox{ and } \left\{x_i = D' \st 1 \leq i \leq n, D' \in \f{1}{2N!}\Z\right\}$$
partition $\R^n$ into (open) standard orthosimplices with edge lengths $\f{1}{2N!}$, see Figure~\ref{fig:partition}.

\mk

We may consider the refinement of Lang's cell decomposition of $E(X)$, obtained by considering all possible hyperplanes $\{d(\cdot,x) \pm d(\cdot,y) = D\}$, for $x,y \in X$ and $D \in \f{1}{N!}\Z$, and $\{d(\cdot,x) = D'\}$, for $x \in X$ and $D' \in \f{1}{2N!}\Z$. Each cell from Lang's decomposition is now refined into a finite union of orthoscheme simplices with edge lengths $\f{1}{2N!}$. Let us denote by $O_NX$ the corresponding simplicial complex. Note that the geometric realization of $O_NX$ is naturally identified with $E(X)$. 

\mk

The vertex set of $O_NX$ will be denoted $X'_N$, and called the Helly subdivision of $X$. When $E(X)$ is realized as an isometric subset of $\R^X$, the vertex set $X'_N$ is naturally identified with
$$X'=E(X) \cap \left(\f{1}{2N!}\N\right)^X= \left\{p \in E(X) \st \forall x \in X, d(p,x) \in \f{1}{2N!}\N\right\}.$$
According to~\cite[Theorem~4.4]{helly_groups}, $X'_N$ is a Helly graph (with edge length $\f{1}{2N!}$).

\mk

Now consider simplices $C,C'$ of $O_NX$, and points $p=\sum_{i=1}^n t_ix_i \in C$ and $p'=\sum_{i'=1}^{n'} t'_{i'}x'_{i'} \in C'$, where $x_1,\dots,x_n$ are the vertices of $C$ and $x'_1,\dots,x'_{n'}$ are the vertices of $C'$. Then we have
\beq d(p,p') &=& \sup_{x \in X} |d(p,x)-d(p',x)| \\
&=& \sup_{x \in X} \sum_{i=1}^n \sum_{i'=1}^{n'} t_i t'_{i'} |d(x_i,x)-d(x'_{i'},x)| \\
&=& \max_{x \in X} \sum_{i=1}^n \sum_{i'=1}^{n'} t_i t'_{i'} |d(x_i,x)-d(x'_{i'},x)|.\eeq
Indeed, since $d(x_i,x) \in \f{1}{2N!}\Z$ for each $1 \leq i \leq n$ and each $x \in X$, and similarly $d(x'_{i'},x) \in \f{1}{2N!}\Z$, we deduce that the supremum is a maximum. 
\ep

\section{Explicit description of the first Helly subdivision} \label{sec:first_subdivision}

We described in Theorem~\ref{thm:orthoscheme_complex}, for each $N \geq 1$, the $N^\text{th}$ subdivision of the injective hull of $X$, which is an orthoscheme simplicial complex, and the $N^\text{th}$ Helly subdivision of the Helly graph itself. When $N=1$, we actually have a very simple and explicit description of these first subdivisions.

\mk

If $X$ is a graph, we say that a clique $\sigma \subset X$ is \emph{round} if it is an intersection of balls of $X$. 

We deduce a very simple and explicit characterization of the the orthoscheme subdivision of the injective hull of a Helly graph.

\bthm \label{thm:simplicial_structure_injective_hull_helly}
Let $X$ denote a Helly graph with finite combinatorial dimension, and let $P_X$ denote the poset of all round cliques of $X$, ordered by inclusion. Then $E(X)$ has a canonical simplicial structure isometric to the $\ell^\infty$ orthoscheme realization of the poset $P_X$ (with edge lengths $\f12$).
\ethm

Before passing to the proof of this result, let us mention this very simple description of the combinatorial dimension of a Helly graph.

\bcor \label{cor:identification_combinatorial_dimension}
Let $X$ denote a Helly graph. Then the combinatorial dimension of $X$ coincides with the length of the longest chain of round cliques of $X$.
\ecor

In particular, if $X$ is uniformly locally finite, or if it has a uniform bound on the size of cliques, then $X$ has finite combinatorial dimension.

\mk

The main technical point in the proof of the theorem is the following lemma.

\blem \label{lem:identification_round_cliques}
Let $X$ denote a Helly graph, let $X'$ denote the first Helly subdivision, and let $P_X$ denote the set of round cliques of $X$. The following map is a bijection:
\beq \sigma:X' & \ra & P_X \\
p & \mapsto & \sigma(p)=\bigcap_{x \in X} B(x,\lceil d(p,x) \rceil).\eeq
\elem

\bp
We will first note that, for each $p \in X$, the subset $\sigma(p) \subset X$ is a non-empty clique. The fact that $\sigma(p) \neq \emptyset$ is a direct consequence of the Helly property since, for any $x,y \in X$, we have $\lceil d(p,x) \rceil + \lceil d(p,y) \rceil \geq d(p,x)+d(p,y) \geq d(x,y)$.

\mk

Note that $p \in E(X)$, $p$ takes values in $\left(\f{1}{2}\N\right)$ and $x \in \sigma(p)$. According to Theorem~\ref{thm:explicit_description_injective_hull}, there exists $y \in X$ such that $d(x,p)+d(p,y)=d(x,y)$. Since $d(x,y) \leq \lceil d(p,y) \rceil$, we deduce that $d(p,x)<1$. In particular, for any $x,y \in \sigma(p)$, we have $d(x,y) \leq d(x,p)+d(p,y) <2$, so $\sigma(p)$ is a clique. Hence $\sigma(p) \in P_X$.

\mk

Conversely, let $A \in P_X$ be a round clique, we will define a map $f=f_A:X \ra \f{1}{2}\N$ as follows. If $A=\{x\}$, then $f_A=d(x,\cdot)$. Otherwise if $x \in A$, let $f(x)=\f{1}{2}$. And if $x \in X \bs A$, let $D \in \N$ denote the minimal distance between $x$ and a point of $A$. If $A \subset B(x,D)$, let $f(x)=D$. If $A \subsetneq B(x,D)$, let $f(x)=D+\f{1}{2}$.

\mk

We claim that $f \in \Delta(X)$. Indeed if $x \in A$ and $y \in X \bs A$, then $f(x)+f(y) \geq d(x,y)$. If $x,y \in X \bs A$, let $D=d(x,A)$ and $D'=d(y,A)$.

If there exists $a \in A$ such that $d(x,a)=D$ and $d(y,a)=D'$, then $f(x)+f(y) \geq D+D' \geq d(x,y)$. 

Otherwise let $a,a' \in A$ such that $d(x,a)=D$ and $d(y,a')=D'$. Then $f(x)+f(y)=D+\f{1}{2}+D'+\f{1}{2} = d(x,a)+d(a,a')+d(a',y) \geq d(x,y)$.

\mk

We claim that $f$ is extremal. For any $x \in A$, let $y \in A \bs \{x\}$: we have $f(x)+f(y)=1=d(x,y)$.

Fix $x \in X \bs A$, and let $D=d(x,A)$. Assume first that $A \subset B(x,D)$. By the Helly property, there exists $z \in X$ that is adjacent to every vertex of $A$, and such that $d(z,x)=D-1$. Since $z \not\in A$, there exists $y \in X$ and $D' \in \N$ such that $A \subset B(y,D')$ and $z \not\in B(y,D')$. As a consequence, $f(x)+f(y)=D+D'=d(x,y)$.

Assume now that $A \not\subset B(x,D)$. Let $a \in A \bs B(x,D)$, we have $f(x)+f(a)=D+\f{1}{2}+\f{1}{2}=D+1=d(x,a)$.

So we have proved that $f$ is extremal. Hence the map $f:P_X \ra X'$ is well-defined.

\mk

We will now prove that $\sigma \circ f=\operatorname{id} : P_X \ra P_X$. Fix $A \in P_X$, we will prove that $\sigma(f_A)=A$. For any $x \in A$, we have $A \subset B(x,\lceil f_A(x) \rceil)$, so $A \subset \sigma(f_A)$. Conversely, let $x \in X \bs A$, and let $D=d(x,A)$. If $A \subset B(x,D)$ then $f_A(x)=D$ and there exists $y \in X \bs A$ and $D' \in \N$ such that $A \subset B(y,D')$ and $x \not\in B(y,D')$. Since $f_A(y)=D'$, we deduce that $x \not\in \sigma(f_A)$. Hence $A=\sigma(f_A)$.

\mk

We will now prove that $f \circ \sigma=\operatorname{id} : X' \ra X'$. Fix $p \in X'$, we will prove that $f=f_{\sigma(p)}=p$. If $p \in X$, then $\sigma(p)=\{p\}$ and $f_{\{p\}}=d(p,\cdot)$. If $p \not\in X$ and $x \in \sigma(p)$, then $d(p,x) < 1$, hence $f(x)=\f{1}{2}=d(p,x)$. If $x \in X \bs \sigma(p)$,  let $D=d(x,\sigma(p))$. For some $y \in \sigma(p)$, we have $d(x,y)=D$, so $|d(x,p)-D| \leq d(y,p)=\f{1}{2}$. Therefore we know that $d(x,p) \in \{D-\f12,D,D+\f12\}$: we want to prove that $d(x,p)=f(x)$.

\mk

Assume first that $\sigma(p) \not\subset B(x,D)$. Let $y,z \in \sigma(p)$ such that $d(x,y)=D$ and $d(x,z)=D+1$. Then $d(x,p) \leq d(x,y)+d(y,p)=D+\f12$, and $d(x,p) \geq d(x,z)-d(z,p)=D+\f12$. Hence $d(x,p)=D+\f12=f(x)$.

\mk

Assume now that $\sigma(p) \subset B(x,D)$.

\mk

We will first prove that, for any $x' \in X \bs \sigma(p)$ adjacent to $\sigma(p)$, we have $d(p,x') = 1$. By contradiction, assume that $d(p,x')=\f12$. Since $x' \not\in \sigma(p)$, there exists $y \in X$ such that $x' \not\in B(y,\lceil d(y,p) \rceil)$. Let $D'$ denote the distance between $y$ and $\sigma(p)$: we deduce that $d(y,p) \leq D'$. Since $d(x',y) \leq d(x',p)+d(p,y) \leq \f12+D'$, we conclude that $d(x',y) \leq D'$, which contradicts the assumption on $y$. 

We will now prove that we have $d(p,x) \geq D$. Let $A \subset X$ denote the set of vertices of $B(x,D)$ adjacent to all vertices of $\sigma(p)$. By the Helly property, we may find $y \in B(x,D-1)$ adjacent to all vertices of $A$. In particular, $y$ is adjacent to all vertices of $\sigma(p)$: since $y \not\in \sigma(p)$, there exists $z \in X$ adjacent to all vertices of $\sigma(p)$ such that $d(y,z)=2$. In particular, $z \not\in A$. So $d(x,z) \geq D+1$. We deduce that $d(x,p) \geq d(x,z) - d(p,z) \geq D+1-1=D$. Hence $d(x,p) \geq D$.

\mk

We will finally prove that $d(p,x) = D$. Since $d(x,p) \in \{D-\f12,D,D+\f12\}$, let us assume by contradiction that $d(x,p)=D+\f12$. According the Helly property, there exists $x' \in X$ adjacent to $\sigma(p)$, such that $d(x,x')=D-1$. Since $d(x,p)=D+\f12$, we have $d(x',p)=\f32$. Let $y \in X$ such that $d(p,x')+d(p,y)=d(x',y)$, and let $D'=d(y,\sigma(p))$: according to the previous case, we know that $d(y,p) \geq D'$. However, if $z \in \sigma(p)$, we have $d(p,y) = d(x',y)-d(p,x') \leq d(x',z)+d(z,y)-\f32 = D'-\f12$, which is a contradiction. Hence $d(p,x)=D$.
\ep

We can now finish the proof of Theorem~\ref{thm:simplicial_structure_injective_hull_helly}:

\bp[of Theorem~\ref{thm:simplicial_structure_injective_hull_helly}]
According to Theorem~\ref{thm:orthoscheme_complex} and Lemma~\ref{lem:identification_round_cliques}, we just have to check that edges in $O_1X$ coincide with edges in the geometric realization of $P_X$.

\mk

If $v,w$ are vertices of $O_1 X$ contained in a common simplex $\sigma$ of $O_1 X$, we will prove that the corresponding round cliques $\sigma,\tau \subset X$ are contained in one another. By contradiction, assume that there exists $x \in \sigma \bs \tau$ and $y \in \tau \bs \sigma$. According to the proof of Lemma~\ref{lem:identification_round_cliques}, we deduce that $d(x,\sigma)=\f12$ and $d(x,\tau)=1$, and similarly $d(y,\sigma)=1$ and $d(y,\tau)=\f12$. Hence $\sigma$ and $\tau$ are separated by the hyperplane $\{p \in E(X) \st d(p,x)-d(p,y)=0\}$ of $O_1 X$. This contradicts the assumption that $v,w$ are adjacent vertices of $O_1 X$.

\mk

Conversely, let us consider two round cliques $\sigma,\tau \subset X$ such that $\sigma \subset \tau$, we will prove that they correspond to adjacent vertices of $O_1 X$. It is sufficient to prove that they are not separated by a hyperplane. Note that one can also check, directly from the definition of $f$ in the proof of Lemma~\ref{lem:identification_round_cliques}, that $d(f_\sigma,d_\tau) \leq \f12$.

Let us fix $x \in X$, $D \in \f12\Z$, since $d(\sigma,\tau)=\f12$, we know that $\sigma$ and $\tau$ are not separated by the hyperplane $\{p \in E(X) \st d(p,x)=D\}$.

Let us fix $x,y \in X$, $\eps=\pm 1$ and $D \in \Z$, and assume by contradiction that  $\sigma$ and $\tau$ are separated by the hyperplane $\{p \in E(X) \st d(p,x)+\eps d(p,y) =D\}$. Since $d(\sigma,\tau)=\f12$, this implies that $d(\sigma,x)+\eps d(\sigma,y) =D \pm \f12$ and $d(\tau,x)+\eps d(\tau,y) =D \mp \f12$. It also implies that $|d(\sigma,x)-d(\tau,x)|=\f12$ and $|d(\sigma,y)-d(\tau,y)|=\f12$. According to the proof of Lemma~\ref{lem:identification_round_cliques}, this implies that there exist $p,q \in \N$ such that, for each $z \in \sigma$, we have $d_X(z,x)=p$ and $d_X(z,y)=q$. Thus $d(\sigma,x)=p$ and $d(\sigma,y)=q$, so $d(\sigma,x)+\eps d(\sigma,y) \neq D \pm \f12$. This is a contradiction.

We conclude that $\sigma$ and $\tau$ are adjacent vertices in $O_1X$.
\ep

\mk

\bthm \label{thm:characterizations_first_helly_subdivision}
Let $X$ denote a Helly graph. The following graphs are naturally isomorphic:
\bit
\item The graph with vertex set
$$E(X) \cap \left( \f12 \N\right)^X,$$
with an edge between $f,g \in E(X)$ if and only if $d_\infty(f,g)=\f12$.
\item The graph with vertex set
$$\{\mbox{round cliques of }X\} = X \cup \{\mbox{non-empty intersections of maximal cliques of }X\},$$
with an edge between $\sigma,\tau \subset X$ if and only if $\sigma \cap \tau \neq \emptyset$ and $\sigma \cup \tau$ is a clique of $X$.
\eit
This graph coincides with the \emph{first Helly subdivision} $X'=X'_1$ of $X$ described in Theorem~\ref{thm:orthoscheme_complex}, and the natural map $X \ra X'$ is a $2$-homothetic embedding.
\ethm

\bp
Let us first prove that the set of round cliques coincide with the union of $X$ and of the intersections of maximal cliques of $X$.

\mk

Let $\sigma \subset X$ denote a round clique, not reduced to a vertex, and let $x \in X \bs \sigma$ adjacent to $\sigma$. Since $\sigma$ is round, there exists a ball $B(y,r)$ containing $\sigma$ but not $x$, with $r \geq 1$. Since $X$ is Helly, there exists $z \in X$ adjacent to $\sigma$ such that $d(y,z)=r-1$. Let $\tau$ denote a maximal clique of $X$ containing $\sigma \cup \{z\}$: we have $x \not\in \tau$. Hence $\sigma$ is an intersection of maximal cliques.

Conversely, any vertex of $X$ is a ball of radius $0$. And if $\sigma$ is an intersection of maximal cliques of $X$, we will prove that $\sigma=\cap_{y \in \sigma} B(y,1)$. Assume that $x \in X \bs \sigma$ is adjacent to $\sigma$. By assumption on $\sigma$, there exists a maximal clique $\tau \supset \sigma$ such that $x \not\in \tau$. Then there exists $z \in \tau$ not adjacent to $x$: we deduce that $x \not\in B(z,1)$, while $\sigma \subset B(z,1)$. Hence $\sigma=\cap_{y \in \sigma} B(y,1)$ is a round clique.

\mk

Together with Lemma~\ref{lem:identification_round_cliques}, this concludes the proof of the equalities. 

\mk

We will now prove that the edges of $X'$ correspond to the given description.

\mk

Let us consider two vertices $\sigma,\tau$ of $X'$ such that $\alpha=\sigma \cap \tau \neq \emptyset$ and $\sigma \cup \tau$ is a clique of $X$. Let us denote $\beta \in X'$ a maximal clique containing $\sigma \cup \tau$. Let $x=\f12{(\alpha+\beta)} \in |P_X|$ denote the midpoint of the edge between $\alpha$ and $\beta$: computing distances in $E(X)$, we have $d(\sigma,x)=\f14$ and $d(\tau,x)=\f14$, hence $d(\sigma,\tau)=\f12$.

Conversely, let us consider two vertices $\sigma,\tau$ of $X'$ such that $d(\sigma,\tau)=\f12$. Let us consider a geodesic $\gamma$ in $|P_X|$ from $\sigma$ to $\tau$: we may assume that $\gamma$ starts by an affine segment inside a (minimal) simplex $S$ of $|P_X|$. This affine segment exits $S$ in the codimension $1$ face $S'$ of $S$ opposite $\sigma$. There are two possibilites now:
\bit
\item If $\sigma$ is not the minimum nor the maximum of the chain corresponding to the simplex $S$, then $d(\sigma,S')=\f14$ (measured in $E(X)$).
\item If $\sigma$ is either the minimum or the maximum of the chain corresponding to the simplex $S$, then $d(\sigma,S')=\f12$ (measured in $E(X)$).
\eit
Since $d(\sigma,\tau)=\f12$, we deduce that we are in the first case, and also there exists a simplex $T$ of $|P_X|$ containing $S' \cup \{\tau\}$. Moreover, let $v_0<v_1<\dots<v_k$ denote the chain in $P_X$ corresponding to the simplex $S'$. We have $v_0 < \sigma,\tau$, hence $\sigma \cap \tau \neq \emptyset$. Moreover $\sigma,\tau < v_k$, hence $\sigma \cup \tau$ is a clique.
\ep

\section{Classification of automorphisms of Helly graphs} \label{sec:classification}

We now turn to the study of automorphisms of Helly graphs, and the proof of the classification Theorem~\ref{mthm:classification_isometries}.

\mk

Fix a Helly graph $X$. An automorphism $g$ of $X$ is called:
\bit
\item \emph{elliptic} if $g$ has bounded orbits in $X$.
\item \emph{hyperbolic} if, for some vertex $x \in X$, the map $n \in \Z \mapsto g^n \cdot x \in X$ is a quasi-isometric embedding.
\item \emph{parabolic} otherwise.
\eit

\mk

Note that there exist parabolic isometries. For instance, let $G$ denote a finitely generated group, with an infinite order element $g \in G$ which is distorted in $G$. Then the action of $g$ by automorphisms on the Helly hull of any Cayley graph of $G$ is parabolic. However, we will see these do not exist if the Helly graph has finite combinatorial dimension.

\mk

We now give several simple equivalent characterizations of elliptic groups of automorphisms.

\bpro \label{pro:characterizations_elliptic}
Let $G$ denote a group of automorphisms of a Helly graph $X$. The following are equivalent:
\ben
\item $G$ stabilizes a round clique in $X$,
\item $G$ stabilizes a vertex of the first Helly subdivision $X'$ of $X$,
\item $G$ fixes a point in the injective hull $E(X)$ of $X$ and
\item $G$ has a bounded orbit in $X$.
\een
Such a group is called an \emph{elliptic} group of automorphisms of $X$.
\epro

\bp\
\bit
\item[$1. \Rightarrow 4.$] If $G$ stabilizes a clique in $X$, it is clear that $G$ has a bounded orbit in $X$.
\item[$4. \Rightarrow 3.$] According to \cite[Proposition~1.2]{lang}, if $G$ has a bounded orbit in $X$, then $G$ has a fixed point in $E(X)$.
\item[$3. \Rightarrow 1.$] Let $p \in E(X)$ denote a point fixed by $G$, and let
$$\phi(p)=\bigcap_{x \in X} B(x,\lceil d(x,p) \rceil).$$
According to the proof of Lemma~\ref{lem:identification_round_cliques}, $\phi(p)$ is a round clique of $X$. Since $p$ is fixed by $G$, we deduce that $\phi(p)$ is stabilized by $G$.
\item[$1. \Leftrightarrow 2.$] Vertices of $X'$ are the round cliques of $X$.
\eit
\ep

In the finite combinatorial dimension case, one can see that such an elliptic group fixes a simplex pointwise.

\blem \label{lem:elliptic_fixed_simplex}
Let $G$ denote an elliptic group of automorphisms of a Helly graph $X$ with finite combinatorial dimension, and assume that $G$ fixes a point $p \in E(X)$ contained in a minimal simplex $C$ of the first subdivision $O_1X$ of $E(X)$. Then $g$ fixes $C$ pointwise.
\elem

\bp
We know that $p$ stabilizes the vertex set of $C$. According to Theorem~\ref{thm:characterizations_first_helly_subdivision}, the vertices of $C$ form a chain of round cliques. Since $g$ preserves the incusion of cliques, we deduce that $g$ fixes $C$ pointwise.
\ep

We deduce the following important classification of automorphisms of Helly graphs.

\bthm \label{thm:dichotomy}
Let $X$ be a Helly graph with finite combinatorial dimension. Then any automorphism of $X$ is either elliptic or hyperbolic. Moreover, any automorphism of $X$ has a non-empty minimal set in $E(X)$.
\ethm

\bp
Let $N-1$ denote the combinatorial dimension of $X$. Fix an automorphism $g$ of $X$. Let $D=\inf_{p \in E(X)} d(g \cdot p,p)$.
Consider any $p \in E(X)$ such that $d(p,g \cdot p) \leq D+\f{1}{2N!}$, and let $C$ denote the minimal simplex of the orthoscheme complex $O_NX$ of $X$ containing $p$. We may assume that the dimension of $C$ is minimal.

Since simplices in $O_NX$ have diameter at most $\f{1}{2N!}$, vertices of $C$ and $g \cdot C$ are at most $D+\f{1}{N!}$ apart.

\mk

Let $x_1,\dots,x_n$ denote the vertices of $C$. Let $\alpha=\f{1}{2N!}$. For each vertex $x \in X$, let us consider the map
\beq f_x:C &\ra& \R\\
q & \mapsto & d(q,x)-d(g \cdot q,x).\eeq
Note that, according to Theorem~\ref{thm:orthoscheme_complex}, the function $f_x$ is affine.

\mk

For each $x \in X$ and $1 \leq i \leq n$, we have $f_x(x_i) \in \alpha \Z$. Moreover, for any $1 \leq j \leq n$, we have
\beq |f_x(x_i)-f_x(x_j)| &\leq& |d(x_i,x) - d(g\cdot x_i,x)-d(x_j,x)+d(g \cdot x_j,x)| \\
&\leq& d(x_i,x_j)+d(g \cdot x_i,g \cdot x_j) \\
&\leq& 2d(x_i,x_j) \leq 2\alpha.\eeq

\mk

Since each $|f_x(x_i)|$ is bounded above by the diameter of $C \cup g \cdot C$, we deduce that there is a finite set ${\mathcal F}=\{f_{y_1},\dots,f_{y_p}\}$ such that, for any vertex $x \in X$, we have $f_x \in {\mathcal F}$.

\mk

For any $q \in C$, we have $d(q,g \cdot q) = \max_{f \in {\mathcal F}} f(q)$.

\mk

Let us assume that $p \in C$ is such that the number of functions $f \in {\mathcal F}$ such that $d(p,g \cdot p)=f(p)$ is maximal. Since the dimension of $C$ is minimal, we deduce that $p$ is in the interior of $C$. Then we deduce that there exist linearly independent functions $f_1,\dots,f_r \in {\mathcal F}$ such that
$$\{q \in C \st \forall 1 \leq i \leq r, f_i(q)=f_i(p)\}=\{p\}.$$

\mk

Since $C=\{t \in (\R_+)^n \st t_1+\dots+t_n=1\}$, let us consider $f:t \in \R^n \mapsto t_1+\dots+t_n$. We deduce that
$$\{t \in \R^n \st f(t)=1,f_1(t)=f_1(p),\dots,f_r(t)=f_r(p)\}=\{p\},$$
where we have chosen arbitrary affine extensions of $f_1,\dots,f_r$ to $\R^n$.

\mk

Now, remark that for each $1 \leq i \leq r$, the function $g_i=f_i-f_i(x_1)f$ has coefficients in $\{-2\alpha,-\alpha,0,\alpha,2\alpha\}$. In particular, $p$ is the unique solution of a linear system of $n \leq N$ equations with linear coefficients in $\{-2\alpha,-\alpha,0,\alpha,2\alpha\}$, and with constant coefficients in $\alpha \Z$.

\mk

According to Lemma~\ref{lem:linear_equations_coefficients_Z}, we deduce that $p = \sum_{i=1}^n t_ix_i$, with each $t_i \in \f{\alpha}{2^{1+2^{N}}}\Z$.

\mk

In particular, the infimum $D$ is realized: in other words, the isometry $g$ of $E(X)$ is semisimple.

\mk

Assume that $D=0$, and let $p \in E(X)$ such that $g \cdot p=p$. According to Proposition~\ref{pro:characterizations_elliptic}, $g$ is elliptic.

\mk

Assume now that $D>0$, and let $p \in E(X)$ such that $d(p,g \cdot p)=D$. According to~\cite[Proposition~3.8]{lang}, $E(X)$ has a conical geodesic bicombing. So according to~\cite[Proposition~4.2]{descombes_lang_flats}, for any $n \geq 1$, we have $\min_{q \in E(X)} d(g^n \cdot q,q)=nD$. In particular, for any $n \in \N$, we have $d(p,g^n \cdot p)=nD$. So the orbit map $n \in \Z \mapsto g^n \cdot p \in E(X)$ is a homothetic embedding: the isometry $g$ is hyperbolic.

\mk

This concludes the proof that any automorphism of $X$ is either elliptic or hyperbolic.
\ep

We deduce the following equivalent characterizations of hyperbolic automorphisms. 

\bpro \label{pro:characterizations_hyperbolic}
Let $g$ denote an automorphism of a Helly graph $X$ with finite combinatorial dimension $N$. The following are equivalent:
\ben[1.]
\item $g$ is hyperbolic, i.e for some vertex $x \in X$, the map $n \in \Z \mapsto g^n \cdot x \in X$ is a quasi-isometric embedding.
\item $g$ has a geodesic axis in the injective hull $E(X)$ of $X$.
\item There exists a vertex $x$ of the Helly subdivision $X'$ of $X$ and integers $1 \leq a \leq 2N$ and $L \in \N \bs \{0\}$ such that $\forall n \in \N, d(x,g^{an} \cdot x)=nL$.
\item There exists a vertex $x$ of the $N^\text{th}$ Helly subdivision of $X'_N$ of $X$ and $L \in \N \bs \{0\}$ such that $\forall n \in \N, d(x,g^n \cdot x)=nL$.
\item $g$ has unbounded orbits in $X$.
\een
\epro

\bp\
\bit
\item[$1. \Rightarrow 2.$] According to Theorem~\ref{thm:dichotomy}, the minimal set of $g$ in $E(X)$ is non-empty. According to~\cite[Proposition~3.8]{lang}, $E(X)$ has a conical, geodesic bicombing. According to~\cite[Proposition~4.2]{descombes_lang_flats}, we deduce that the isometry $g$ has a geodesic axis in $E(X)$.
\item[$2. \Rightarrow 3.$]  Let $D=\min_{p \in E(X)} d(g \cdot p,p) > 0$, and let $p \in E(X)$ such that $d(p,g \cdot p)=D$. Since $g$ has a geodesic axis in $E(X)$, we may assume that $p$ lies in a simplex $C$ of $O_1X$ of codimension at least $1$. Let $x_1,\dots,x_n$ denote the vertices of $C$: we have $n \leq (N+1)-1=N$.


\mk

Let us first assume that $D \geq 1$, we will show that $D$ is rational, and its denominator is a divisor of $2(k-1)$, with $k \leq n+1 \leq N+1$.

\mk

 For each $k \geq 2$, let $A_k=\{1,2,\dots,n\}^k$. For each $a \in A_k$, and for every vertex $y \in X$, let us define
$$f_y(a) = \sum_{i=1}^{k-1} |d(g^{i-1} \cdot x_{a_i},y)-d(g^i \cdot x_{a_{i+1}},y)|.$$
This quantity should roughly be thought as the length of a path going through vertices of $C,g \cdot C, \dots g^{k-1} \cdot C$. Let us also define
$$\alpha_y = \inf\left\{\f{f_y(a)}{k-1} \st k \geq 2, a \in A_k, a_1=a_k\right\}.$$
More precisely, one can interpret these quantities in terms of maximal lengths of paths in a graph as follows. Consider the finite graph $\Gamma$ with vertices labeled $1,\dots,n$, such that given any two vertices $i,j$, there exists one oriented edge from $i$ to $j$, whose length depend on a time parameter $t \geq 1$: at time $t$, its length is $|d(x_i,g^{-t+1} \cdot y)-d(g \cdot x_j,g^{-t+1} \cdot y)|$. The set $A_k$ is the set of oriented paths of $k$ vertices in $\Gamma$, with time $1 \leq t \leq k-1$, and $f_y(a)$ is the length of the path $a$. Finally, $\alpha_y$ is the minimal average length of an oriented loop.

\mk

We claim that $\alpha_y$ is attained by some element $a \in A_k$ with $a_1=a_k$ such that $k \leq n+1$. Consider some $k \geq 3$ and $a \in A_k$ with $a_1=a_k$ such that, for any $k' < k$ and $a' \in A_{k'}$ with $a'_1=a'_{k'}$, we have $\f{f_y(a')}{k'-1} > \f{f_y(a)}{k-1}$. We will prove that $k \leq n+1$. By contradiction, if $k > n+1$, since there are $n$ vertices $x_1,\dots,x_n$, there exists a strict subloop $a'$ of $a$ consisting of $k'$ vertices, with $k'<k$. Since $\f{f_y(a')}{k'-1} > \f{f_y(a)}{k-1}$, removing the loop $a'$ decreases the average length of the loop, which contradicts the assumption. Hence $k \leq n+1$.

\mk

Any two vertices of $X'$ have distance in $\f{1}{2}\N$, and since $A_{n+1}$ is finite, we also deduce that $\alpha_y$ is attained, and furthermore $\alpha_y \in \f{1}{2(k-1)}\N$, for some $k \leq n+1$. In particular $\alpha_y \in \f{1}{2N!}\N$.

\mk

Let $j \geq 2$ and $a \in A_j$ with $a_1=a_j$ such that $\f{f_y(a)}{j-1} = \alpha_y$. Without loss of generality, we may assume that $j$ is large enough such that, if there exists $q \in \f{1}{2N!}\N$ such that $|D-q| \leq \f{1}{j-1}$, then $D=q$. We can also assume that $j-1$ is a multiple of $2N!$. According to Theorem~\ref{thm:explicit_description_injective_hull}, there exists $y \in X$ such that $njD=d(p,g^{nj} \cdot p)=d(p,y)-d(g^{nj} \cdot p,y)$.

\mk

Note that
\beq njD = d(p,g^{nj} \cdot p) & = &  d(p,y)-d(g^{nj} \cdot p,y) \\
& = & \sum_{i=1}^{nj} d(g^{i-1} \cdot p,y)-d(g^i \cdot p,y).\eeq
Since $|d(g^{i-1} \cdot p,y)-d(g^i \cdot p,y)| \leq d(g^{i-1} \cdot p,g^i \cdot p)=D$, for any $1 \leq i \leq nj$, we have $d(g^{i-1} \cdot p,y)-d(g^i \cdot p,y)=D$.

\mk

Moreover, for any $1 \leq i \leq nj$ and $(a_1,a_2) \in A_2$, we have
\beq d(x_{a_1},g^{1-i} \cdot y) - d(g \cdot x_{a_2},g^{1-i} \cdot y) & \geq & d(p,g^{1-i} \cdot y) - d(g \cdot p,g^{1-i} \cdot y)-d(x_{a_1},p)-d(g \cdot x_{a_2},g \cdot p) \\
& \geq & d(p,g^{1-i} \cdot y) - d(g \cdot p,g^{1-i} \cdot y)-1 \\
& \geq & d(g^{i-1} \cdot p,y) - d(g^i \cdot p,y)-1 \\
& \geq & D-1 \geq 0,\eeq
since $D \geq 1$ by assumption.

\mk

Also remark that, since $a_j=a_1$, we have
$$|d(x_{a_1},y)-d(g^{j-1} \cdot x_{a_1},y)| \leq \sum_{i=1}^{j-1} |d(g^{i-1} \cdot x_{a_i},y)-d(g^i \cdot x_{a_{i+1}},y)|  \leq f_y(a) = (j-1)\alpha_y.$$
Note that $|d(x_{a_1},y)-d(g^{j-1} \cdot x_{a_1},y)| \geq |d(p,y)-d(g^{j-1} \cdot p,y)| - 1 = (j-1)D - 1$. Hence $(j-1)D -1 \leq (j-1)\alpha_y$, and so $D \leq \alpha_y + \f{1}{j-1}$.

\mk

Consider the integer $h=nj+1$. For any $a \in A_h$, there exists a subloop consisting of at least $h-n$ vertices, hence $f_y(a) \geq (h-n)\alpha_y$. According to Theorem~\ref{thm:orthoscheme_complex}, let $t_1,\dots,t_n \in \R_+$ such that $t_1+\dots+t_n=1$ and $p=t_1x_1+\dots+t_nx_n$. We have
\beq &&\sum_{a \in A_h} t_{a_1}t_{a_2} \cdots t_{a_h} f_y(a) \\
&=& \sum_{a \in A_h} t_{a_1}t_{a_2} \cdots t_{a_h} (|d(x_{a_1},y)-d(g \cdot x_{a_2},y)| + |d(g \cdot x_{a_2},y)-d(g^2 \cdot x_{a_3},y)| + \\
& & \dots + |d(g^{h-2} \cdot x_{a_{h-1}},y)-d(g^{h-1} \cdot x_{a_h},y)|) \\
&=& \sum_{i=1}^{h-1}  \sum_{a \in A_h} t_{a_1}t_{a_2} \cdots t_{a_h} |d(g^{i-1} \cdot x_{a_1},y)-d(g^i \cdot x_{a_2},y)| \\
&=& \sum_{i=1}^{h-1} \sum_{a \in A_2} t_{a_1}t_{a_2} |d(x_{a_1},g^{1-i} \cdot y)-d(g \cdot x_{a_2},g^{1-i} \cdot y)| \\
&=& \sum_{i=1}^{h-1} \sum_{a \in A_2} t_{a_1}t_{a_2} (d(x_{a_1},g^{1-i} \cdot y)-d(g \cdot x_{a_2},g^{1-i} \cdot y)) \\
&=& \sum_{i=1}^{h-1}\left( \sum_{\ell=1}^n t_\ell d(x_\ell,g^{1-i} \cdot y)-\sum_{\ell=1}^n t_\ell d(g \cdot x_\ell,g^{1-i} \cdot y)\right)\\
&=& \sum_{i=1}^{h-1} (d(p,g^{1-i} \cdot y)-d(g \cdot p,g^{1-i} \cdot y))\\
&=& d(p,y) - d(g^{h-1} \cdot p,y)=(h-1)D.\eeq

For any $a \in A_h$, we have $f_y(a) \geq (h-n)\alpha_y$, hence $(h-1)D \geq (h-n)\alpha_y$. So $D \geq \alpha_y - \f{n-1}{h-1}$. We deduce that
$$|D-\alpha_y| \leq \max\left(\f{n-1}{h-1},\f{1}{j-1}\right) = \max\left(\f{n-1}{nj},\f{1}{j-1}\right) \leq \f{1}{j-1},$$
so by choice of $j$ we conclude that $D=\alpha_y$.

\mk

In particular, $D$ is rational, and its denominator is a divisor of $2(k-1)$, with $k \leq n+1 \leq N+1$.

\mk

Assume now that $D$ is not necessarily greater than $1$. Choose arbitrary distinct prime integers $q,q' \geq N+1$ such that $qD,q'D \geq 1$. According to the previous argument applied to $g^q$ and $g^{q'}$, we deduce that $D$ is rational, and its denominator is a divisor of $2q(k-1)$ and of $2q'(k'-1)$, for some $k,k' \leq n+1 \leq N+1$. Hence we conclude that the denominator 
of $D$ is a divisor of $2(k-1)$, with $k \leq n+1 \leq N+1$.

\mk

Now let us consider a vertex $z \in X$ and $1 \leq i_0 \leq n$ such that $d(x_{i_0},z)-d(g^{k-1} \cdot x_{i_0},z)$ is maximal.

If there exist $1 \leq i,i' \leq n$ such that $d(x_i,z)-d(g^{k-1} \cdot x_{i'},z) \leq (k-1)D-\f12$, then $d(x_i,z)-d(g^{k-1} \cdot x_{i},z) \leq (k-1)D$. By maximality of $z$, we deduce that for every $z' \in X$ we have $d(x_{i},z')-d(g^{k-1} \cdot x_{i},z') \leq (k-1)D$, hence $d(x_i,g^{k-1} \cdot x_i)=(k-1)D$. In particular, $x_i$ lies on an axis for $g^{k-1}$.

Otherwise, for all $1 \leq i,i' \leq n$, we have $d(x_i,z)-d(g^{k-1} \cdot x_{i'},z) \geq (k-1)D$ since vertices in $X'$ have distances in $\f12\N$. Now
$$(k-1)D = d(p,g^{k-1} \cdot p) \geq d(p,z)-d(g^{k-1} \cdot p,z) = \sum_{i,i'=1}^n t_i t_{i'} \left(d(x_i,z)-d(g^{k-1} \cdot x_{i'},z)\right),$$
so we conclude that, for all $1 \leq i,i' \leq n$, we have 
$$d(x_i,z)-d(g^{k-1} \cdot x_{i'},z) = (k-1)D.$$
In particular, we have $d(x_1,g^{k-1} \cdot x_1)=(k-1)D$, so $x_1$ lies on an axis for $g^{k-1}$.

Either way, there exists $1 \leq i \leq n$ such that the vertex $x_i \in X'$ satisfies that for any $r \in \N$, we have $d(x_i,g^{2(k-1)r} \cdot x_i) = r2(k-1)D$, with $L=2(k-1)D \in \N$.

\item[$2. \Rightarrow 4.$] Using the notations from the proof of $2. \Rightarrow 3.$, let us assume that $p$ lies in a simplex $C$ of $O_1X$ of minimal dimension denoted $n-1$, with vertices $x_1,\dots,x_n$ in $X'$. According to $3.$, we know that $D \in \f{1}{2n!}\N$. Let us furthermore assume that $p$ is as close as possible from a vertex of $X'_n$: let us be more precise.

\mk

According to~\cite[Theorem~4.5]{lang}, there exist vertices $z_1,\dots,z_{n-1}$ of $X$ such that the map $q \in C \mapsto (d(q,z_1),\dots,d(q,z_{n-1})) \in (\R^{n-1},\ell^\infty)$ is an isometric embedding. Moreover, given any $q \in C$ and $q' \in EX$, there exists $1 \leq j \leq n-1$ such that $d(q,q')=|d(q,z_j)-d(q',z_j)|$. Let us now assume precisely that $p \in C$ is such that the number of $1 \leq j \leq n-1$ for which $d(p,z_j) \in \f{1}{2n!}\N$ is maximal. We will prove that in fact $p$ is a vertex of $X'_n$.

Let us assume, without loss of generality, that $z_1,\dots,z_r$ are such that, for all $1 \leq j \leq r$, we have $d(p,z_j)-d(g \cdot p,z_j)=D$ and, for all $r+1 \leq j \leq n-1$, we have $|d(p,z_j)-d(g \cdot p,z_j)|<D$. For each $1 \leq j \leq r$, let $1 \leq j' \leq n-1$ such that $g^{-1} \cdot z_j$ is equivalent to $z_{j'}$ with respect to $C$, i.e. there exists $\eps_j=\pm1$ and $a_j \in \f{1}{2}\Z$ such that, for any $q \in C$, we have $d(q,g^{-1} \cdot z_j)=\eps_j d(q,z_{j'})+a_j$. We deduce that, for all $1 \leq j \leq r$, we have $d(p,z_j)-\eps_j d(p,z_{j'})=D \in \f{1}{2n!}\N$.

If $r<n-1$, we may find $p' \in C$ with a larger number of coordinates in $\f{1}{2n!}\Z$. Hence we deduce that $r=n-1$, and so for all $1 \leq j \leq r$ we have $d(p,z_j) \in \f{1}{2n!}\N$. So $p \in X'_n$. In particular, $p$ is a vertex of the $N^\text{th}$ Helly subdivision.
\item[$3. \Rightarrow 5.$] This is immediate.
\item[$4. \Rightarrow 5.$] This is immediate.
\item[$5. \Rightarrow 1.$] If $g$ has unbounded orbits in $X$, by definition $g$ is not elliptic. According to Theorem~\ref{thm:dichotomy}, $g$ is hyperbolic.
\eit
\ep

We deduce the following interesting corollary about translations lengths in Helly graphs, which directly generalizes the analogous theorem by Gromov about translation lengths in Gromov-hyperbolic groups (see~\cite[8.5.S]{gromov_hyperbolic_groups}). Since Garside groups are Helly according to~\cite{huang_osajda_helly}, this implies a direct analogue of~\cite{lee_lee_garside_translation} for a very closely related translation length.

\bcor
Let $X$ denote a Helly graph with finite combinatorial dimension $N$. Then any hyperbolic automorphism of $X$ has rational translation length in $X$, with denominator uniformly bounded by $2N$.
\ecor

\section{Fixed points for pairs of elliptic subgroups} \label{sec:fixed_point_pair}

We now use the orthoscheme subdivision complexs to study fixed point sets of pairs of elliptic subgroups, that is used in~\cite{haettel_osajda_locally_elliptic} for the study of locally elliptic actions on Helly graphs.

\bpro \label{pro:fixed_GH}
Let $X$ denote a Helly graph with finite combinatorial dimension $N-1$, and let $G,H$ denote elliptic automorphism groups of $X$. Then the distance between the fixed point sets $E(X)^G$ and $E(X)^H$ is realized by vertices in the Helly subdivision $X'_{2N}$ of $X$.
\epro

\bp
The proof will be very similar to that of Theorem~\ref{thm:dichotomy}. Let $p \in E(X)^G$ and $p' \in E(X)^H$. Denote by $C,C'$ the minimal simplices of $O_N X$ containing $p,p'$ respectively. Without loss of generality, assume that the dimensions of $C$ and $C'$ are minimal. According to Lemma~\ref{lem:elliptic_fixed_simplex}, we know that $C \subset E(X)^G$ and $C' \subset E(X)^H$. Let $\alpha=\f{1}{2N!}$. Let us denote the vertices of $C$ (resp. $C'$) by $x_1,\dots,x_n$ (resp. $x'_1,\dots,x'_{n'}$).

\mk

For each vertex $x \in X$, let us consider the map
\beq f_x:C \times C' &\ra& \R\\
(q,q') & \mapsto & d(q,x)-d(q',x).\eeq
Note that, according to Theorem~\ref{thm:orthoscheme_complex}, the function $f_x$ is affine.

\mk

For each $x \in X$, $1 \leq i \leq n$ and $1 \leq i' \leq n'$, we have $f_x(x_i,x'_{i'}) \in \alpha \Z$. Moreover, for any $1 \leq j \leq n$, we have
$$|f_x(x_i,x'_{i'})-f_x(x_j,x'_{i'})| \leq |d(x_i,x)-d(x_j,x)| \leq d(x_i,x_j) \leq \alpha.$$
Similarly, for $1 \leq j' \leq n'$, we have $|f_x(x_i,x'_{i'})-f_x(x_i,x'_{j'})| \leq \alpha$.

\mk

We deduce that there is a finite set ${\mathcal F}=\{f_{y_1},\dots,f_{y_p}\}$ such that, for any vertex $x \in X$, we have $f_x \in {\mathcal F}$.

\mk

For any $(q,q') \in C \times C'$, we have $d(q,q') = \max_{f \in {\mathcal F}} f(q,q')$.

\mk

Let us assume that $(p,p') \in C \times C'$ is such that the number of functions $f \in {\mathcal F}$ such that $d(p,p')=f(p,p')$ is maximal. Since the dimensions of $C$ and $C'$ are minimal, we deduce that $(p,p')$ is in the interior of $C \times C'$. Then we know that there exist linearly independent functions $f_1,\dots,f_r \in {\mathcal F}$ such that
$$\{(q,q') \in C \times C' \st \forall 1 \leq i \leq r, f_i(q,q')=f_i(p,p')\}=\{(p,p')\}.$$

\mk

Since $C=\{t \in (\R_+)^n \st t_1+\dots+t_n=1\}$ and $C'=\{t' \in (\R_+)^{n'} \st t'_1+\dots+t'_{n'}=1\}$, let us consider $f:(t,t') \in \R^n \times \R^{n'} \mapsto t_1+\dots+t_n$ and $f':(t,t') \in \R^n \times \R^{n'} \mapsto t'_1+\dots+t'_{n'}$. We deduce that
$$\{(t,t') \in \R^n \times \R^{n'} \st f(t,t')=1,f'(t,t')=1,f_1(t,t')=f_1(p,p'),\dots,f_r(t,t')=f_r(p,p')\}=\{(p,p')\},$$
where we have chosen any affine extension of $f_1,\dots,f_r$ to $\R^n \times \R^{n'}$.

\mk

Now, remark that for each $1 \leq i \leq r$, the function $g_i=f_i-f_i(x_1,x'_1)(f+f')$ has coefficients in $\{-\alpha,0,\alpha\}$. In particular, $(p,p')$ is the unique solution of a linear system of $n+n' \leq 2N$ equations with linear coefficients in $\{-\alpha,0,\alpha\}$, and with constant coefficients in $\alpha \Z$.

\mk

According to Lemma~\ref{lem:linear_equations_coefficients_Z}, there exists $1 \leq D \leq N!$ such that $p = \sum_{i=1}^n t_ix_i$, with each $t_i \in \f{\alpha}{D}\Z$, and similarly  $p' = \sum_{i=1}^{n'} t'_ix'_i$, with each $t'_i \in \f{\alpha}{D}\Z$.

\mk

We deduce that the distance $d(E(X)^G,E(X)^H)$ is attained, and it is realized by vertices of $X'_{N^2}$, since $\f{D}{\alpha}=2N!D$ divides $2N!^2$, which itself divides $2(2N)!$.\ep

\blem \label{lem:linear_equations_coefficients_Z}
Let us consider a matrix $A \in GL(n,\mathbb{Q})$, such that each coefficient of $A$ is in $\{-1,0,1\}$, and let $y \in \Z^n$. Then $A^{-1}y \in \left(\f{1}{D}\Z\right)^n$, where $D \geq 1$ divides $n!$.
\elem

\bp
The determinant of $A$ is such that $D=|\det(A)| \leq n!$. Therefore each coefficient of $A^{-1}y$ lies in $\f1D\Z$.
\ep

\bibliography{../../../bibli}
\bibliographystyle{alpha}
	
\end{document}